\documentclass[12pt]{amsart}
\usepackage{graphicx}
\usepackage[T1]{fontenc}
\usepackage[english,francais]{babel}

\usepackage{amssymb,url,xspace}
\usepackage{amsfonts}
\usepackage{amsmath}
\usepackage{latexsym}
\usepackage{amscd}
\usepackage{color}
\usepackage{shadow}
\usepackage{a4wide}
\usepackage{endnotes}
\usepackage{amsopn}
\usepackage{url}
\usepackage{dsfont}
\usepackage{pdfsync}

\theoremstyle{plain}
\newtheorem{thm}{Théorème}
\newtheorem{lemma}{Lemme}
\newtheorem{rmq}[thm]{Remarque}
\newtheorem{prop}{Proposition}

\title[Estimations de Strichartz]{Estimations de Strichartz pour l'équation des ondes dans un domaine strictement  convexe}
  \author{Oana Ivanovici}
  \address{Laboratoire J. A. Dieudonn\'e, UMR CNRS 7351\\
    CNRS et Universit\'e Nice C\^ote d'Azur\\
    Parc Valrose\\
    06108 Nice Cedex 02\\
    France} \email{oana.ivanovici@unice.fr}
    
  \author{Gilles Lebeau}
  \address{Laboratoire J. A. Dieudonn\'e, UMR CNRS 7351\\
    Universit\'e Nice C\^ote d'Azur\\
    Parc Valrose\\
    06108 Nice Cedex 02\\
    France} \email{gilles.lebeau@unice.fr}

  \author{Fabrice Planchon}
  \address{Laboratoire J. A. Dieudonn\'e, UMR CNRS 7351\\
    Universit\'e Nice C\^ote d'Azur\\
    Parc Valrose\\
    06108 Nice Cedex 02\\
    France}
  \email{fabrice.planchon@unice.fr} \thanks{ authors were
    partially supported by A.N.R. grant GEODISP and ERC project SCAPDE} 

\begin{document}

\begin{abstract}
On se propose d'établir ici des estimations de Strichartz pour l'équation des ondes dans un domaine strictement convexe (quelconque) de $\mathbb{R}^3$. 
Dans le papier \cite{ilpdisp}, nous avons obtenu des estimations de dispersions optimales pour la solution de l'\'equation des ondes dans un domaine convexe particulier (le mod\`ele de Friedlander). Ce r\'esultat, qui montre qu'une perte de $\frac 14$ par rapport au cas plat est n\'ecessaire, implique, en utilisant la méthode $TT^*$ usuelle, des estimations de type Strichartz sans perte (en termes d'échelle) mais avec des indices modifiés, ce que l'on résume en parlant de perte d'un quart. Pour réussir à faire mieux (obtenir un résultat qui correspondrait à une perte d'au plus $\frac 16$), il faut s'intéresser au lieu et à la fréquence de l'apparition des caustiques responsables de la perte et montrer qu'elles sont suffisamment exceptionnelles pour que l'effet d'une moyenne en temps, présent dans les estimations de Strichartz, puisse atténuer la perte.
\end{abstract}


\maketitle
\tableofcontents 
\vfill\break

\section{Introduction}

Les estimations dispersives dites "de Strichartz" mesurent la taille et
la dispersion des solutions de l'équation des ondes linéaire sur
un domaine $\Omega$ avec bord $\partial\Omega$ (possiblement vide):
\begin{equation} \label{WE} \left\{ \begin{array}{l} (\partial^2_t-
\Delta) u(t, x)=0, \;\; x\in \Omega \\ u|_{t=0} = u_0 \; \partial_t
u|_{t=0}=u_1,\\ Bu=0,\quad x\in \partial\Omega.
 \end{array} \right.
 \end{equation} Ici $\Delta$ désigne l'opérateur de
Laplace-Beltrami sur $\Omega$. Si $\partial\Omega\neq \emptyset$, on
considère ou bien la condition de Dirichlet sur le bord ($B = $
l'opérateur identité): $u|_{\partial\Omega}=0$ ou bien la
condition de Neumann ($B=\partial_{\nu}$), où $\nu$ est le vecteur
unitaire normal au bord.

Pour pouvoir perturber ces équations et étudier les problèmes
non-linéaires associés, avoir un contrôle de la
"taille" du flot linéaire en termes de la taille des données
initiales s'avère crucial. Pour l'équation des ondes
non-linéaires, les normes mixtes $L^q_tL^r(\Omega)$ sont
particulièrement utiles: au prix d'une moyenne en temps on gagne de
l'intégrabilité en espace, parfois jusqu'à $r=\infty$.

Quant aux équations linéaires, une estimation (locale) de base
indique que sur toute variété riemannienne sans bord, la solution
linéaire de \eqref{WE} vérifie (pour $T<\infty$)
\begin{equation}\label{stricrd} h^{\beta}\|\chi(hD_t)u\|_{L^q([0,T],
L^r(\mathbb{R}^d))}\leq C\Big(\|u(0,x)\|_{L^2}+\|hD_t u\|_{L^2}\Big),
\end{equation} où $\chi\in C^{\infty}_0$ est une fonction lisse à
support inclus dans un voisinage de $1$. Si $d$ désigne la dimension
de la variété, on a $\beta=d\Big(\frac 12-\frac 1r\Big)-\frac 1q$,
où le couple $(q,r)$ est admissible pour l'équation des ondes,
c.a.d.:
\begin{equation}\label{adm} \frac 2q+\frac{d-1}{r}\leq
\frac{d-1}{2},\quad  q>2.
\end{equation} Lorsque l'égalité a lieu dans \eqref{adm}, la paire
$(q,r)$ est dite strictement admissible. Si \eqref{stricrd} a lieu
pour $T=\infty$, on parle d'inégalité de Strichartz globale en
temps. Ces estimations ont été étudiées depuis bien longtemps dans
l'espace de Minkowski (métrique plate):  si $\Omega$
désigne l'espace $\mathbb{R}^d$ avec la métrique euclidienne
$g_{i,j}=\delta_{i,j}$, la solution $u_{\mathbb{R}^d}(t,x)$ de
\eqref{WE} dans $\mathbb{R}^d$ avec $(u_0=\delta_a, u_1=0)$ est
donnée par la formule
 \[ u_{\mathbb{R}^d}(t,x)=\frac{1}{(2\pi)^d}\int
\cos(t|\xi|)e^{i(x-a)\xi}d\xi
 \] et elle vérifie les estimations de dispersion usuelles:
\begin{equation}\label{disprd}
\|\chi(hD_t)u_{\mathbb{R}^d}(t,.)\|_{L^{\infty}(\mathbb{R}^d)}\leq
C(d) h^{-d}\min\{1, (h/t)^{\frac{d-1}{2}}\}.
\end{equation} L'interpolation entre \eqref{disprd} et l'estimation de
l'énergie, suivie d'un argument classique de dualité dit $TT^*$, permet
d'obtenir facilement les estimations \eqref{stricrd}. Ces estimations
peuvent être généralisées à tout $(\Omega,g)$ sans bord
grâce à leur caractère local (vitesse de propagation finie). Les
estimations \eqref{stricrd} sont optimales sur une telle variété
riemannienne.

La motivation principale pour les estimations de type Strichartz vient
de leurs applications en analyse harmonique et l'étude des
problèmes non-linéaires dispersifs. Par exemple, \eqref{stricrd}
peut être utilisée pour montrer des résultats d'existence pour
l'équation des ondes semi-linéaire.

Même si le cas sans bord est relativement bien compris depuis un certain
temps, l'obtention de tels résultats sur des variétés à bord
s'avère une tache bien plus difficile. Pour des variétés à
bord strictement concave, cette théorie a pu être établie
grâce à la paramétrice de Melrose et Taylor près de rayons
tangents au bord: des estimations de Strichartz optimales pour les ondes à l'extérieur d'un obstacle strictement convexe ont été obtenues dans \cite{smso95} et, très récemment, des estimations de dispersion optimales en $d=3$ ont été établies dans \cite{il16}, ainsi que des contre-exemples en dimension plus grande à l'extérieur d'une sphère. 
Pourtant, dès que l'hypothèse de stricte
concavité du bord est enlevée, la présence des rayons
géodésiques multi-réfléchis et de leurs limites, les rayons
glissants, ne permet plus d'avoir une telle paramétrice. En dehors du
cas d'un bord concave, il n'y avait que très peu de résultats jusqu'à très
récemment: des estimations de Strichartz avec pertes ont été obtenues dans \cite{blsmso08} dans un domaine compact, en utilisant les constructions de paramétrices en temps petit de \cite{smso07}, qui, à leurs tour, ont été inspirées des travaux sur des domaines à métriques à régularité faible \cite{tat02}. L'avantage majeur de \cite{blsmso08} est en même temps son point faible: en considérant uniquement des intervals de temps qui ne permettent pas de voir plus d'une réflection d'un paquet d'onde au bord, on peut traiter le problème indépendamment de la géométrie du bord mais on ne peut pas voir l'effect de dispersion dans les directions tangentes (en dimension $d\geq 3$).  

Dans cette note qui résume les idées importantes  de \cite{illp} et \cite{ilpstri}, le but est d'obtenir des
estimations de Strichartz à l'intérieur d'un domaine, meilleures que
celles obtenues directement à partir du résultat (optimal) de dispersion
de \cite{illp} :
avant d'énoncer notre résultat principal, on va introduire le
modèle de Friedlander d'un demi-espace $\Omega_d=\{(x,y)| x>0,
y\in\mathbb{R}^{d-1}\}$ muni de la métrique $g_F$ héritée de
l'opérateur de Laplace suivant:
$\Delta_F=\partial^2_x+(1+x)\Delta_{\mathbb{R}^{d-1}_y}$. On s'aper\c
coit facilement que $(\Omega_d,g_F)$ modélise localement un domaine
strictement convexe: en effet, $(\Omega_d,g_F)$ peut être regardé
comme un modèle simplifié du disque unité $D(0,1)$ après le
passage en coordonnées polaires $(r,\theta)$, avec $r=1-x/2$,
$\theta=y$. Pour ce mod\`ele isotrope particulier, nous avons montré dans \cite{ilpdisp} qu'une perte de
dérivées par rapport à l'estimation de dispersion libre
\eqref{disprd} est inévitable, et qu'elle appara\^it en raison de la
présence de caustiques de type queue d'aronde dans le support
singulier de la solution; ensuite, nous avons obtenu dans \cite{ilpstri} des estimations de Strichartz avec perte de moins de $\frac 16$.

\begin{thm}\label{thmdisper}\cite{ilpdisp} Il existe $T>0$ et il
existe une constante $C(d)>0$ tels que pour tous $a\in (0,1]$, $h\in
(0,1]$ et $t\in (0,T]$ la solution $u_a(t,x,y)=\cos
(t\sqrt{|\Delta_F|})(\delta_{x=a,y=0})$ de \eqref{WE} avec $\Delta=\Delta_F$ vérifie
\begin{equation}\label{dispco} |\chi(hD_t)u(t,x)|\leq C(d)
h^{-d}\min\{1,(h/t)^{\frac{d-2}{2}}\gamma(t,h,a)\},
\end{equation} où
 \[ \gamma(t,h,a)=\left\{\begin{array}{c} (\frac ht)^{1/2}+
a^{1/4}(\frac ht)^{1/4}$, si $a\geq h^{4/7-\epsilon}\\ (\frac
ht)^{1/3}+h^{1/4}$, si $a\leq h^{1/2}.
 \end{array} \right.
 \] De plus, il existe une suite de temps $t_n=4n\sqrt{a(1+a)}$
pour lesquels on a égalité dans \eqref{dispco} pour $x=a$.
\end{thm}
\begin{rmq} L'estimation \eqref{dispco} nous dit que dans un domaine
strictement convexe on perd une puissance $\frac 14$ dans l'exposant de $h$ par rapport
à l'estimation \eqref{stricrd} de l'espace libre, ce qui est due à
des phénomènes micro-locaux comme les caustiques générées en
temps arbitrairement petit près du bord. Ces caustiques
apparaissent lorsque les rayons optiques envoyés d'une même
source dans des directions différentes cessent de diverger.
\end{rmq}

Notre but est de généraliser les résultats précédents au cas d'un
domaine strictement convexe quelconque. Pour cela, on refait d'abord la
construction de la paramétrice en suivant la méthode de notre papier
\cite{ilpstri}, mais cette fois dans le cadre d'un opérateur de
Laplace général. Cette étape implique plusieurs difficultés
techniques importantes : l'idée ``simple'' consistant à dire qu'on
sera proche du cas modèle anisotrope ne trouve pas de traduction
élégante dans une preuve qui permettrait de tordre un opérateur sur
l'autre. Il faut donc effectuer une construction microlocale
perturbative ``à la main''.

\begin{thm} \label{thmstriilp}\cite{illp} Les estimations de
Strichartz restent vraies pour la solution de \eqref{WE} dans un domaine strictement convexe $\Omega\subset \mathbb{R}^3$ avec 
\[
\frac 1q=\Big(\frac{d-1}{2}-\frac
16\Big)\Big(\frac 12-\frac 1r\Big), \quad d=3.
\] 
\end{thm}
\begin{rmq} Ce résultat a été démontré en dimension $d=2$
par M.Blair, H.Smith et C.Sogge dans \cite{blsmso08} pour des
métriques arbitraires (i.e. sans l'hypothèse de convexité
stricte). Notre théorème de \cite{ilpstri} améliore tous les résultats connus
jusqu'à présent pour $d\geq 3$. Le cas d'un strict convexe
quelconque, esquissé ici, est traité dans un travail en cours en collaboration avec
G.Lebeau, F.Planchon et R.Lascar \cite{illp}.
\end{rmq}

Un point essentiel de la preuve consiste en une description précise
de la géométrie des ondes "sphériques" (la "sphère" étant
ici un objet singulier en raison des multiples réflexions). La
démonstration fournit une analyse assez détaillée de la fonction
de Green des ondes, au moins dans certains régimes. En effet, la
construction des paramétrices microlocales que j'ai utilisée pour
obtenir les contre-exemples de \cite{doicef}, \cite{doi2} semble instable dans une zone
intermédiaire près du bord, où l'estimation de dispersion est
obtenue par injection de Sobolev, en remarquant que le flot préserve
essentiellement la taille du microsupport.  

\section{Paramétrice pour un domaine modèle}

Dans cette partie on va indiquer comment construire une paramétrice dans le cas d'un opérateur modèle anisotrope. Soit $\Omega_d$ définit plus haut. On introduit l'opérateur de Laplace suivant
\[
\Delta_M=\partial^2_x+\sum\partial^2_{y_j}+x\Big(\sum_{j,k}r_{j,k}\partial_{y_j}\partial_{y_k}\Big),
\]
avec la condition de Dirichlet sur le bord. On note $q(\eta)=\sum_{j,k} r_{j,k}\eta_j\eta_k$, où $r_{j,k}$ sont tels que $q$ est définie positive. Sous cette hypothèse $(\Omega_d,\Delta_M)$ modélise localement un domaine strictement convexe qui coïncide avec le modèle de Friedlander lorsque $q(\eta)=|\eta|^2$. En prenant la
transformation de Fourier dans la variable transverse $y$, $-\Delta_M$
devient $-\partial^2_x+\eta^2+xq(\eta)$, qui, pour $\eta\neq 0$ est
auto-adjoint et positif sur $L^2(\mathbb{R}_+)$ avec résolvante
compacte. Il admet une base orthonormale dans $L^2(\mathbb{R}_+)$ de
fonctions propres $\{e_k(x,\eta)\}_{k\geq 0}$ associées aux valeurs
propres $\lambda_k(\eta)=\eta^2+\omega_kq(\eta)^{2/3}$, où
$\{-\omega_k\}_{k\geq 1}$ désignent les zéros de la fonction
d'Airy en ordre décroissant. On a une formule explicite
\begin{equation} e_k(x,\eta)=f_k\frac{q(\eta)^{1/6}}{k^{1/6}}
Ai\Big(q(\eta)^{1/3}x-\omega_k\Big),
\end{equation} où, pour $k\geq 1$, $f_k$ est
tel que $\|e_k(.,\eta)\|_{L^2(\mathbb{R}_+)}=1$, 
$ \int_0^{\infty}Ai^2(x-\omega_k)dx
=\frac{k^{1/3}}{|f_k|^2}$.

Pour $a>0$, soit $\delta_{x=a}$ la distribution de
Dirac sur $\mathbb{R}_+$: alors elle s'écrit comme une somme de
modes $e_k$ de la fa\c con suivante
\[ \delta_{x=a}=\sum_{k\geq 1} e_k(x,\eta)e_k(a,\eta).
\] 
On considère la donnée au temps $t=0$ de la forme
$u_0(x,y)=\psi(hD_y)\delta_{x=a,y=0}$, où $h\in (0,1]$ est un petit
paramètre et où $\psi\in C^{\infty}_0((\frac 12,2))$.

On rappelle la formule de la fonction de Green associée au temps $t$
obtenue dans \cite{ilpdisp} à partir de la décomposition de la
donnée $u_0$ en somme de modes de galerie:
\begin{equation}\label{greenfct} G_M(x,y,t,a)=\sum_{k\geq
1}\frac 1h\int_{\mathbb{R}}e^{\pm it\sqrt{\lambda_k(\eta/h)}}
e^{iy\eta/h} 
 \psi(\eta)e_k(x,\eta/h)e_k(a,\eta/h)d\eta.
\end{equation}
\begin{rmq} La fonction $G_M$ est une paramétrice pour l'équation
$i\partial_t u\pm\sqrt{-\Delta_M}u=0$ (et donc de l'équation des ondes \eqref{WE}) avec donnée
au temps $t=0$ égale à $u_0$. Elle est valable pour toute distance
initiale $a>0$ au bord.
\end{rmq}

\begin{rmq}
Notons que dans la somme \eqref{greenfct}, la contribution principale vient des valeurs $k\simeq \frac{a^{3/2}}{h}$. Pour $k\ll \frac{a^{3/2}}{h}$ ou $k\gg \frac{a^{3/2}}{h}$ il est facile d'obtenir des estimations de dispersion avec perte de $\frac 16$. Les modes $\frac 1h \lesssim k$ correspondent à des ondes transverses pour lesquelles il n'y a pas de perte dans la dispersion.
\end{rmq}
\subsubsection*{L'estimation de dispersion} On considère deux cas, selon la taille de $a$.
\\

\begin{itemize}
\item Pour de petites valeurs de la distance initiale au bord $a$,
$h^{2/3}\lesssim a\ll h^{22/39}$, on utilise directement la
formule explicite de $G_M$ pour estimer la solution de
l'équation des ondes dans $\Omega$ en norme $L^{\infty}$. Dans
\cite[Théorème 1.8.(2)]{doicef}, on a démontré que si la
donnée est un mode de galerie, i.e. de la forme
\begin{equation}\label{mdgk} u^k_0(x,y)=\frac 1h \int
e^{iy\eta/h}\psi(\eta)e_k(x,\eta/h)e_k(a,\eta/h)d\eta,
\end{equation} avec $k\geq 1$ fixé, alors on n'obtient pas de perte
dans les estimations de Strichartz par rapport à l'espace
$\mathbb{R}^d$ avec métrique plate (mais la constante dépend bien sûr
de $k$). On utilise la même
méthode pour montrer que si $a$ est assez petit, étant donné la somme
définissant $G_{M}$ il n'y a pas "trop" de termes, on peut obtenir
la dispersion avec une perte d'au plus $1/6$.

\item Pour $a\gg h^{4/7}$ on obtient une paramétrice sous la forme
d'une somme indexée selon le nombre de réflexions au
bord. Chaque terme dans la somme est une intégrale dont la phase
admet des points critiques dégénérés et il se trouve que si
$a>h^{1/3}$, la perte correspondante dans l'estimation $L^{\infty}$
est bien plus importante pour qu'on puisse espérer appliquer un
argument $TT^*$ classique pour obtenir des Strichartz optimales. Dans
ce régime, on localise la fonction de Green près de points
où une singularité de type queue d'aronde se forme et on obtient
des estimations plus raffinées dans de très petits voisinages
autour de ces points.

\begin{rmq} Les deux régimes se recouvrent bien car $h^{4/7}\ll
h^{22/39}$.
 \end{rmq}

\begin{thm}\label{thmbon} \cite{ilpstri} Dans le régime $a\gg h^{4/7}$, le facteur $h^{1/4}$ du Théorème \ref{thmdisper} apparaît seulement près d'une suite de points $t_n$, avec une
estimation (optimale) de $\gamma(t,h,a)$ pour $t$ dans l'intervalle
$I_n=(t_n(1-a), t_n(1+a))$ :
\begin{equation}\label{newgamma} \gamma(t,h,a)\leq (\frac
ht)^{1/2}+h^{1/3}+\frac{a^{1/8}h^{1/4}}{n^{1/4}+h^{-1/12}a^{-1/24}|t^2-t_n^2|^{1/6}}.
\end{equation} Notons aussi que pour $t\notin I_n$, le dernier facteur
est $\leq (h/t)^{1/3}$. Ce raffinement de $\gamma(t,h,a)$ est donné par
une analyse soignée de la dégénérescence des arguments de phase
stationnaire autour de $t_n$ dans \cite{ilpdisp}.
\end{thm}
\end{itemize}

\section{Paramétrice dans le cas d'un domaine général}

Soit $\Omega\subset\mathbb{R}^d$, $d\geq 2$ un domaine strictement convexe \`a bord $C^{\infty}$, et $\Delta$ le Laplacien dans $\Omega$ avec condition de Dirichlet sur $\partial\Omega$. Dans un syst\`eme de  coordonn\'ees g\'eod\'esiques normales $(x,y)\in\Omega_d=\mathbb{R}_+\times\mathbb{R}^{d-1}$, la m\'etrique de $\mathbb{R}^d$ est \[
\left(
\begin{array}{cc}
g(x,.)  &   0   \\
 0 &     1
\end{array}
\right),
\]
o\`u $g$ est la m\'etrique Riemannienne induite sur l'hypersurface dist$((x,y),\partial\Omega_d)=x$. Dans ces coordonn\'ees le Laplacien s'\'ecrit sous la forme
\[
\Delta=\partial_{x}^2+ R(x,y,\partial_{y}), \quad (x,y)\in\mathbb{R}_+\times \mathbb{R}^{d-1},
\] 
avec $R_0(y,\partial_y)=R(0,y,\partial_y)=\sum\partial_{y_j}\partial_{y_k}+O(y^2)$ et $R_{1}(y,\partial_{y})= \partial_{x}R(0,y,\partial_{y})=\sum R_{1}^{j,k}(y)
\partial_{y_{j}}\partial_{y_{k}}$ (où, grâce à la condition de convexité stricte on sait que la forme quadratique $q(\eta)=\sum R_{1}^{j,k}(y)\eta_{j}\eta_{k}$ est définie positive). 

Pour $a>0$, on cherche \`a obtenir une solution approchée de 
\begin{equation}\label{ondeeq}
\partial^2_t u-\Delta u=0,\quad (u,\partial_t u)|_{t=0}=(\delta_{x=a,y=0},0),\quad u|_{\mathbb{R}\times\partial\Omega_d}=0,
\end{equation}
 à partir d'une solution (explicite!) de l'opérateur modèle donné par: 
\[
\triangle_{M}= \partial_{x}^2+ \sum \partial_{y_{j}}^2+ x\Big( \sum R_{1}^{j,k}(0)
\partial_{y_{j}}\partial_{y_{k}}\Big).
\]
Pour $\Delta=\Delta_M$, une solution de \eqref{ondeeq} peut \^etre facilement obtenue de fa\c con explicite. Le passage de la solution de l'équation des ondes avec Laplacien modèle (anisotrope) $\Delta_M$ au cas général de $\Delta$ nécessite l'utilisation du théorème des surfaces {\it glancing} de Richard Melrose:
\begin{thm}
(Melrose-Taylor, Eskin) Deux paires de hypersurfaces avec intersection "glancing" sont localement équivalentes, via une transformation canonique symplectique $\chi_M$.
\end{thm}
Dans le cas de l'opérateur modèle $\Delta_M$,  les hypersurfaces $\{q_M=x=0\}$ et $\{p_M=\xi^2+\eta^2+xq(\eta)-1=0\}$ ont une intersection {\it glancing} au point $(x,y,\xi,\eta)=(0,0,0,1)$, c.a.d. qu'elles vérifient la condition suivante:
\begin{equation}\label{glancingcond}
\{p,q\}=0, \{p,\{p,q\}\}\neq 0\text{ and }\{q,\{p,q\}\}\neq 0\text{ en } (0,0,0,1).
\end{equation}
Dans le cas qui nous concerne, la paire 
\[
\{q=X=0\}\text{ and }\{p=\Xi^2+R(X,Y,\Theta)-1=0\}
\]
vérifie aussi la condition \eqref{glancingcond}. Le Théorème de Melrose nous dit qu'il existe une transformation canonique $\chi_M$ qui vérifie 
\[
\chi_M(x=0, \xi^2+\eta^2+xq(\eta)=1)=(X=0, \Xi^2+R(X,Y,\Theta)=1).
\]
A partir d'une fonction génératrice de $\chi_M$ on va pouvoir obtenir une paramétrice du cas général à partir d'une paramétrice du cas modèle. En fait, on peut trouver explicitement des symboles $p_{0,1}$ tels que, si on introduit
\begin{equation*}
G_M(x,y,\eta,\omega)=e^{iy\eta}\Big(p_0Ai(xq(\eta)^{1/3}-\omega)+xp_{1}\vert\eta\vert^{-1/3}Ai'(xq(\eta)^{1/3}-\omega)\Big),
\end{equation*}
alors $G_M$ vérifie $-\triangle_M G_M= \rho^2 G_M + O_{C^\infty}(\vert\eta\vert^{-\infty})$.

\subsubsection*{Fonction g\'en\'eratrice} Soit $(X-x)u+(Y-y)v+\Gamma(X,Y,u,v)$ une fonction génératrice de $\chi_M$, alors on peut montrer qu'il existe un symbole $p(x,y,\eta,\omega,\sigma)$ de degré $0$, à support compact dans un voisinage de $(0,0,0,\eta)$, $\eta\in\mathbb{R}^{d-1}\setminus 0$ pour lequel
\[
G(x,y;\eta,\omega)= {1\over 2\pi}\int e^{i(y\eta+\sigma^3/3+\sigma(xq(\eta)^{1/3}-\omega)
+\rho\Gamma(x,y, {\sigma q(\eta)^{1/3}\over\rho},{\eta\over\rho}))}p\ d\sigma
\]
vérifie $-\Delta G= \rho^2 G + O_{C^\infty}(\vert\eta\vert^{-\infty})$. A partir de $G$ et à l'aide d'une formule sommatoire de type Airy-Poisson on va obtenir une solution de \eqref{ondeeq} entre deux réflexions successives qui va avoir des propriétés similaires (points critiques dégénérés du même ordre) que celle du cas modèle. Si on note $\Phi$ la phase de $G$ et $p_h$ une normalisation de $p$ alors
 \[
K_\omega(f)(t,x,y)={h^{2/3}\over 2\pi}\int e^{{i\over h}(t\rho(h^{2/3}\omega,\theta)+\Phi-y'\theta-t'h^{2/3}\omega)}p_{h}f(y',t') \ dy'dt'd\theta ds
\]
est une solution de \eqref{ondeeq} entre $t=0$ et la première réflexion au bord et il existe $f_a$ tel que $K_\omega(f_a)(0,x,y)=\delta_{x=a,y=0}$. Le but est d'obtenir des paramétrices de la même forme entre deux réflexions successives; dans ce but on utilise une formule de type Poisson.

\subsubsection*{Formule d'Airy-Poisson} On pose $A_{\pm}(z)=e^{\mp i
\pi/3}Ai(e^{\mp i\pi/3}z)$. On introduit la fonction suivante de $\mathbb{R}$
dans $\mathbb{C}$ :
\[ L(\omega)=\pi+i\log {\Big(\frac{A_-(\omega)}{A_+(\omega)}\Big)}.
\]
\begin{lemma}\label{AiPoi} La fonction $L$ est à valeurs réelles, analytique, strictement croissante
et vérifie:
\[ L(0)=\pi/3,\quad \lim_{\omega\rightarrow -\infty}L(\omega)=0, \quad
L(\Omega) \sim_{\omega\rightarrow +\infty}\frac 43\omega^{3/2},
\] et pour tout $k\geq 1$ on a
\[ L(\omega_k)=2\pi k\Leftrightarrow Ai(-\omega_k)=0, \quad
L'(\omega_k)=\int_0^{\infty}Ai^2 (x-\omega_k)dx.
\]
Ici $\{-\omega_k\}_{k\in\mathbb{N}}$ d\'esignent les zeros de la fonction d'Airy $Ai$ en ordre d\'ecroissant.
\end{lemma}
Le lemme se démontre en utilisant les expansions asymptotiques
associées à $A_\pm$ et des calculs élémentaires.
\begin{prop} On a l'égalité suivante dans l'espace des
distributions $\mathcal{D}'(\mathbb{R}_{\omega})$:
\begin{equation}\label{egaldistrib} \sum_{N\in\mathbb{Z}}
e^{-iNL(\omega)}=2\pi \sum_{k\in\mathbb{N}^*}\frac{1}{L'(\omega_k)}
\delta_{\omega=\omega_k}.
\end{equation}
\end{prop}
\begin{rmq} Notons que l'égalité précédente résulte simplement de la
  classique formule de Poisson reliant la somme des $\exp(-iN x)$ au
  peigne de Dirac, suivi d'un changement de variable $x=L(\omega)$,
  qui permet d'indexer la somme de droite sur les zéros d'Airy et non pas sur les
nombres naturels (ce qui va rendre nos calculs ultérieurs beaucoup
plus limpides).
\end{rmq}

On pose maintenant:
\begin{align}\label{parametrix}
\mathcal P_{h,a}(t,x,y) & = <\sum_{N\in \mathbb Z} e^{-iNL(\omega)},
K_{\omega}(g_{h,a})>_{\mathcal D'(\mathbb R)}\\
\nonumber
 & = 2\pi \sum_{k\in \mathbb N^*}{1\over L'(\omega_{k})}K_{\omega_{k}}(g_{h,a}).
\end{align}

\begin{prop}
$\mathcal P_{h,a}(t,x,y)$ est une paramétrice de \eqref{ondeeq} qui vérifie la condition de Dirichlet.
\end{prop}

On a ainsi obtenu une solution sous la forme d'une somme d'intégrales oscillantes (la somme sur $N$) qui sont presque-orthogonales en temps. Estimer la norme $L^{\infty}$ de la somme à $t$ fixé revient à estimer le $\sup$ des normes $L^{\infty}$ de chaque terme. Chaque terme est une intégrale oscillante avec un unique point critique d'ordre $3$ qui apparaît uniquement pour $x=a$ et une suite $(t_n,y_n)$ qui n'est plus explicite mais dépend des directions initiales. On obtient alors l'équivalent du Théorème \ref{thmbon} dans ce cas général \cite{illp}.

Notons que les termes de la somme sur $k\in \mathbb{N}$ dans \eqref{parametrix} sont les "modes de galerie".  Ces modes n'apparaissent dans la littérature que dans le cas d'un domaine modèle (la boule unit\'e ou un mod\`ele de Friedlander). Dans le travail \cite{illp} on construit des modes de galerie (\`a partir de la formule de $K_{\omega_k}(g_{h,a})$) dans le cas d'un op\'erateur g\'en\'eral pour $k$ suffisamment grand, ce qui nous permet d'agir comme dans \cite{ilpdisp} et de d\'emontrer des estimations de dispersion de fa\c con directe pour des valeurs tr\`es petites de la distance initiale au bord $a$.

\section{Les estimations de Strichartz optimales dans le régime des
queues d'arondes} Pour simplifier on se restreint au cas de la
dimension $d=3$ et $\Delta=\Delta_F$.  On considère la fonction de Green
$G(t,x,y,a)=\chi(hD_t) e^{it\sqrt{|\Delta_F|}}(\delta_{x=a,y=0})$ et
pour $f$ à support compact dans les variables $(s,a\geq 0,b)$, on
pose
\[ A(f)(t,x,y)=\int G(t-s,x,y-b,a)f(s,a,b)ds da db.
\] L'exposant dispersif est dans ce cas
$\alpha_d:=\frac{d-1}{2}-\frac{1}{6}=\frac 56$. Il s'agit d'estimer la
norme de Strichartz $L^{12/5}([0,1],L^{\infty}(\Omega_3))$ (donc
$r=\infty$ et $q=12/5$):
\[ h^{2\beta}\|A(f)\|_{L^{12/5}_{t\in[0,1]}L^{\infty}_{x,y}}\leq
C\|f\|_{L^{12/7}_s L^1_{a,b}},\quad 2\beta=(d-\alpha_d)=3-5/6=13/6.
\] On résume la situation: les singularités de type queue d'aronde
apparaissent seulement en $t_n$, $x=a$; elles ont un effet sur les intervalles de
temps $I_n:=(t_n(1-a),t_n(1+a))$.  En dehors de $I_n$ on ne voit que
des cusps qui font perdre $(\frac{h}{t})^{-1/6}$ dans la dispersion et
induisent donc les estimations de Strichartz avec $q=12/5$. L'estimation de $\gamma(t,h,a)$ dans \eqref{newgamma}
permet de se localiser précisément là où l'argument usuel de
type $TT^*$ ne s'applique plus.

On écrit $G(t,x,y,a)=G_0(t,x,y,a)+G_s(t,x,y,a)$ où $G_s$ dénote
la partie singulière, associée à une localisation en espace -
temps de $G$ dans des boules centrées aux points où les queues
d'arondes apparaissent, i.e. en
\[ |x-a|\leq \frac{a}{n^2},\quad |t-t_n|\leq
a^{3/2}{n}.
\]
 En utilisant la section précédente, on obtient les estimations raffinées suivantes:
\begin{prop}\label{propraf}
\[ h^{2\beta}\sup_{x,y}|G_0(t,x,y,a)|\leq C|t|^{-5/6};
\]
\[ h^{2\beta}\sup_{x,y}|G_s(t,x,y,a)|\leq D(t,a,h),\quad
\sup_{a,h}\int_{-1}^1|D(t,a,h)|^pdt<\infty,\quad \forall p<3.
\]
\end{prop} Soit $A=A_0+A_s$, le découpage correspondant à la
décomposition précédente. L'estimation pour $A_0$ en découle
facilement, car la convolution par $|t|^{-5/6}$ envoie $L^{12/7}$ dans
$L^{12/5}$. En utilisant la Proposition \ref{propraf} on déduit que
$h^{2\beta}A_s$ est borné de $L^1_s L^1_{a,b}$ dans
$L^{3-\epsilon}_t L^{\infty}_{x,y}$: remarquons qu'il est
indispensable de faire la convolution ($L^p\ast L^1 \rightarrow L^p$)
avant d'intégrer en $a$ :
\begin{align*}
  \| A_s(f)(t,x,y) \|_{L^\infty_{x,y}}  \leq & \int
  \sup_{x,y}|G_s(t-s,x,y-b,a)|\times |f(s,a,b)|dsdadb\\
 \lesssim &  \int
  h^{-2\beta} D(t-s,a,h)| |f(s,a,b)| ds da db\\
  \| A_s(f)(t,x,y) \|_{L^{p}_t L^\infty_{x,y}}  \lesssim &  h^{-2\beta}\int  \| D(t-s,a,h)\|_{L^p_t} \| f(s,a,b)\|_{L^1_s} da db\\
  \| A_s(f)(t,x,y) \|_{L^{p}_t L^\infty_{x,y}}  \lesssim &  h^{-2\beta} \sup_{a,h} \| D(t,a,h)\|_{L^p_t} \| f(s,a,b)\|_{L^1_{s,a,b} }\,.
\end{align*}

Comme on travaille
avec des normes locales en temps, on en déduit immédiatement
l'estimation souhaitée puisque $1<12/7$ et $12/5<3$.

\end{document}